\documentclass[12pt,leqno]{article} 
\usepackage{amsmath,amssymb,latexsym,theorem,epsfig}
\usepackage[all]{xy}

\makeatletter
\let\@internalcite\cite
\def\cite{\def\citeauthoryear##1##2{##1{##2}}\@internalcite}
\def\@biblabel#1{\def\citeauthoryear##1##2{##1{##2}}[#1]\hfill}
\makeatother

\theoremstyle{plain} 

\newtheorem{thm}{Theorem}
\newtheorem{theo}{Theorem}[section]
\newtheorem{lem}[theo]{Lemma}
\newtheorem{corr}[theo]{Corollary}

\newtheorem{defi}[theo]{Definition}
\newtheorem{prop}[theo]{Proposition}
\newtheorem{prop-defi}[theo]{Proposition-Definition}
\newtheorem{lemma-defi}[theo]{Lemma-Definition}
{\theorembodyfont{\rmfamily} \newtheorem{rem}[theo]{Remark}}
{\theorembodyfont{\rmfamily} }

\newcommand{\op}[1]{\operatorname{#1}}

\newcommand{\liesp}{{\mathfrak s}{\mathfrak p}}
\newcommand{\ang}{\boldsymbol{d}\boldsymbol{a}}

\begin{document}

\centerline{\LARGE\sf Hyperelliptic Szpiro inequality}

\

\bigskip
\bigskip

\noindent
{\large  Fedor Bogomolov}\footnote{Partially supported by NSF Grant
DMS-9801591}   \\
Courant Institute of Mathematical sciences \\
New York University \\
New York, NY 10012

\medskip

\noindent
{\large Ludmil Katzarkov}\footnote{Partially supported by NSF Career
Award DMS-9875383 and A.P. Sloan research fellowship} \\
University of California at Irvine \\ 
Irvine, CA 92697

\medskip

\noindent
{\large Tony Pantev}\footnote{Partially supported by 
NSF Grant DMS-9800790 and A.P. Sloan research fellowship} \\
University of  Pennsylvania \\ 
209 South 33rd street \\ 
Philadelphia, PA 19104-6395

\

\bigskip

\

\begin{abstract} We generalize the classical Szpiro inequality to the
case of a semistable family of hyperelliptic curves. We show that for
a semistable symplectic Lefschetz fibration of hyperelliptic curves
of genus $g$, the number $N$ of non-separating vanishing cycles and
the number $D$ of singular fibers satisfy the inequality $N \leq (4g+2)D$.
\end{abstract}

\tableofcontents

\addtolength{\baselineskip}{4pt}

\section{Introduction} \label{sec-intro}

The classical Szpiro inequality \cite{SZ} asserts that for any semistable
algebraic family of genus one curves $f : X \to {\mathbb C}{\mathbb
P}^{1}$, the number of components of the singular fibers of $f$ is
bounded from above by $6$ times the number of singular fibers. A
symplectic generalization of Szpiro's result was proven in \cite{abkp} by a
purely group-theoretic technique. Unfortunately the
analogous bounds for fibrations of higher genus curves are extremely
hard to obtain (or even guess) already in the algebraic-geometric
setup. 

In this note we generalize the techniques developed in \cite{abkp}
to obtain a proof of a Szpiro type bound for symplectic families of
hyperelliptic curves. For hyperelliptic curves over number fields such
a bound was conjectured by P.Lockhart in \cite{lockhart}. Our goal is
to prove the following symplectic version of Lockhart's conjecture.

\begin{thm} \label{thm-main}
Let $f : X \to S^{2}$ be a symplectic fibration of hyperelliptic
curves of genus $g$  with only semi-stable fibers. Assume further that
$f$ admits a topological section and that all the
vanishing cycles of $f$ are non-separating. Let $D$ be the number of
singular fibers of $f$ and let $N$ be the number of vanishing cycles.
Then $N \leq (4g+2) D$.
\end{thm}

Note that as a special case of this theorem one obtains a Szpiro
inequality for algebraic families of hyperelliptic curves over
${\mathbb C}{\mathbb P}^{1}$. From a slightly different perspective
the Szpiro inequality can be viewed as an
obstruction for the existence of a symplectic structure on a the total
space of a topological Lefschetz fibration. Indeed,
Theorem~\ref{thm-main} implies that a topological Lefschetz fibration
which violates the inequality $N \leq (4g+2) D$ can not be
symplectic or, equivalently, orientable (see e.g. \cite{abkp},
\cite{gompf-stipsicz}). 

The paper is organized as follows. In section~\ref{sec-setup} we
recall some (mostly standard) material about hyperelliptic  symplectic
fibrations, the hyperelliptic mapping class group and its relation
with the braid group. In section~\ref{sec-Abg} we describe an
criterion for the triviality of a central extension of an
Artin braid group. This criterion plays an important role in the proof
of Theorem~\ref{thm-main} - it provides an efficient way of
controlling the ambiguity in lifting relations from the hyperelliptic
mapping class group to the braid group. In section~\ref{sec-angles} we
introduce our main technical tool - the displacement angle of an
element in the universal cover of the symplectic group. Finally in
section~\ref{sec-proof} we compare the values of the degree character
and the displacement angle character on the braid group and use this
comparison to deduce the hyperelliptic Szpiro inequality. We conclude
in section~\ref{sec-remarks} with a brief discussion of some
ideas concerning the general Szpiro inequality.

\bigskip

\noindent
{\bf Acknowledgments.} We would like to thank N.Ivanov for hepful
discussions on the subject of this paper, A.Agboola for bringing
Lockhart's paper to our attention and S.-W. Zhang for explaining to us
his proof of the elliptic Szpiro inequality for pencils over bases of
higher genus.

\section{Hyperelliptic symplectic fibrations} \label{sec-setup}

First we recall some basic definitions and results and describe the
precise setting in which the Szpiro inequality will be considered in
this paper. More details can be found in the papers \cite{abkp},
\cite{gompf-stipsicz}, \cite{siebert-tian}, \cite{ivan-thesis}.

Let $(X,\omega)$ be a smooth compact symplectic 4-fold. A {\em differentiable
fibration} on $X$ is a surjective $C^{\infty}$ map $f : X \to S^{2}$ with
finitely many critical points $Q_{1}, Q_{2}, \ldots, Q_{N}$ (not
necessarily in distinct fibers) such that locally near each $Q_{i} \in
X$ and $f(Q_{i}) \in S^{2}$,
there exist complex analytic coordinates $x,y$ on $X$ and $t$ on
$S^{2}$, so that $t = f(x,y) = x^{2} + y^{2}$. A differentiable
fibration $f : X \to S^{2}$ is called {\em symplectic} if the
smooth fibers of $f$ are symplectic submanifolds with respect to
$\omega$ and if for every $Q_{i}$ the symplectic form $\omega_{Q_{i}}
\in \wedge^{2}T_{Q_{i}}^{*}X$ is non-degenerate on each of the two
planes contained in the tangent cone of $f^{-1}(f(Q_{i}))$ at $Q_{i}$. 
In particular, for a symplectic fibration the local complex analytic
coordinates around each $Q_{i}$ can be chosen to be compatible with a
global orientation on $X$.

For a point $p \in S^{2}$ we will denote the fiber $f^{-1}(p)$ by
$X_{p}$. Since by definition the rank of $df$ drops only at the points
$Q_{i}$, it follows that for each $p$ the fiber $X_{p}$ is singular
only at the points $X_{p}\cap \{Q_{1}, \ldots, Q_{N} \}$. Let $p \in
S^{2}$ and let $X_{p}^{\sharp} = X_{p}-\{Q_{1}, \ldots, Q_{N}\}$ be
the smooth locus of $X_{p}$. A compact surface $Z \subset X$ which is
the closure of some connected component of some $X_{p}^{\sharp}$ is
called a {\em fiber component} of $f : X \to S^{2}$. Note that for
each $p$ the homology class $[X_{p}] \in H_{2}(X,{\mathbb Z})$ splits
as a sum $[X_{p}] = \sum_{\Sigma \in \pi_{0}(X_{p}^{\sharp})}
n_{\Sigma}[\overline{\Sigma}]$, where $\overline{\Sigma}$ denotes the
closure of $\Sigma$ 
in $X$ and $n_{\Sigma}$ is a positive integer - the {\em multiplicity} of
the fiber component $\overline{\Sigma}$. Again the assumption that the
$Q_{i}$'s are the only critical points of $f$ implies that $n_{\Sigma}
= 1$ for all possible fiber components. 

Let $f : X \to S^{2}$ be a symplectic fibration of fiber genus $g
\geq 1$. By analogy with the algebro-geometric case we will say that
$f$ is {\em semistable} if and only if for every $p \in S^{2}$ and
every $\Sigma \in \pi_{0}(X_{p}^{\sharp})$ of genus zero we have that
$\Sigma$ is homeomorphic to a sphere with at least two punctures.

Given a symplectic fibration $f : X \to S^{2}$, we denote by $p_{1},
\ldots, p_{D} \in S^{2}$ the critical values of $f$. The restriction
of $f$ to $S^{2}-\{p_{1},\ldots,p_{D}\}$ is a $C^{\infty}$ fiber
bundle with a fiber some closed oriented surface $C_{g}$ of genus
$g$. Choose a base point $o \in S^{2}-\{p_{1},\ldots,p_{D}\}$ and put 
$\op{mon} : \pi_{1}(S^{2}-\{p_{1},\ldots,p_{D}\},o) \to \op{Map}_{g}
:= \pi_{0}(\op{Diff}^{+}(X_{o}))$ for the
corresponding geometric monodromy representation.

The hyperelliptic fibrations are singled out among all possible
symplectic fibrations by a condition on the geometric monodromy. Fix
a double cover $\nu : C_{g} \to S^{2}$ and let $\iota \in
\op{Map}_{g}$ denote the mapping class of the covering involution. The
{\em hyperelliptic mapping class group} of genus $g$ is the centralizer
$\Delta_{g}$ of $\iota$ in $\op{Map}_{g}$:
\[
\Delta_{g} := \{ \phi \in \op{Map}_{g} | \phi\iota\phi^{-1} = \iota \}. 
\]
Similarly we can consider versions of $\Delta_{g}$ that take into
account punctures on $C_{g}$. Concretely, denote by
$\op{Diff}^{+}(C_{g})_{r}^{n}$  the group of
all orientation preserving diffeomorphisms of $C_{g}$ preserving $n+r$
distinct points on $C_{g}$ and inducing the identity on the tangent
spaces at $r$ of those points. Let $\op{Map}_{g,r}^{n} :=
\pi_{0}(\op{Diff}^{+}(C_{g})_{r}^{n})$ and define
\[
\Delta_{g,r}^{n} := \Delta_{g}\times_{\op{Map}_{g}} \op{Map}_{g,r}^{n}.
\]
With this notation we can now define
\begin{defi} \label{defi-hesf} A {\em hyperelliptic symplectic
fibration} on a smooth symplectic $4$-fold $(X,\omega)$ is a
symplectic fibration $f : X \to S^{2}$, with a monodromy
representation is conjugate to a representation taking values in
$\Delta_{g}$. 

The fibration $f: X \to S^{2}$ is said to be a {\em hyperelliptic
symplectic fibration with a section}, if $f$ has a topological section
and the corresponding monodromy representation 
in $\op{Map}_{g}^{1}$
is conjugate to one taking values in $\Delta_{g}^{1}$.
\end{defi}

A classical theorem of Kas \cite[Theorem~2.4]{kas} asserts that for a
symplectic fibration $f : X \to S^{2}$ of genus $g \geq 2$ the
diffeomorphism type of $f$ is uniquely determined by the geometric
monodromy of $f$. Therefore the fact that $f : X \to S^{2}$ is a hyperelliptic
fibration with a section is equivalent to the existence of a
topological section of $f$ together with an involution of $X$ which
preserves the section and acts as a hyperelliptic involution on each
fiber of $f$.

In the remainder of this paper we will consider only semistable
hyperelliptic fibrations with a section. The geometric monodromy
representation for such an $f : X \to S^{2}$ sends a small closed
loop running once counterclockwise 
around one of the $p_{i}$ into the product of right
handed Dehn twists about the 
cycles vanishing at the points $\{Q_{1},\ldots, Q_{N}\} \cap X_{s_{i}}$. 
Thus the monodromy representation $\op{mon} : \pi_{1}(S^{2}-\{p_{1},\ldots,
p_{D}\}) \to \Delta_{g}$ is
encoded completely in the relation in $\Delta_{g}$: 
\[
\tau_{1}\tau_{2}\ldots \tau_{N} = 1,
\]
where $\tau_{i} \in \Delta_{g}$ denotes the mapping class of the 
right-handed Dehn twist in $\op{Diff}^{+}(C_{g})$
about the loop vanishing at $Q_{i}$. 

Similarly the monodromy
representation $\pi_{1}(S^{2}-\{p_{1},\ldots,
p_{D}\}) \to \Delta_{g}^{1}$ is completely encoded in the relation in
the group $\Delta_{g}^{1}$:
\[
{\mathfrak t}_{1}{\mathfrak t}_{2}\ldots {\mathfrak t}_{N} = 1
\]
where ${\mathfrak t}_{i} \in \Delta_{g}^{1}$ denotes the right-handed Dehn twist
$\op{Diff}^{+}(C_{g})^{1}$ about the loop vanishing at $Q_{i}$.

\section{Central extensions of Artin braid groups} \label{sec-Abg}

In this section we recall some standard facts about Artin braid groups
and study an important class of central extensions of such groups.

Let $\Gamma$ be a graph with a vertex set $I$. Assume that $\Gamma$
has no loops and that any two vertices of $\Gamma$ are connected by at
most finitely many edges.

\begin{defi} \label{defi-Abg}
The {\em Artin braid group} associated with $\Gamma$ is the group
$\op{Art}_{\Gamma}$  generated by elements $\{ t_{i} | i \in I
\}$, so that if $i, j \in I$ are two distinct vertices connected by $k_{ij}$
edges, then $t_{i}$ and $t_{j}$ satisfy the relation
\[
t_{i}t_{j}t_{i}t_{j}\ldots = t_{j}t_{i}t_{j}t_{i}\ldots,
\]
where both sides are words of length $k_{ij} + 2$.
\end{defi}

\

\begin{rem} \label{rem-Abg} (i) Note that by specifying a graph
$\Gamma$ one specifies not only the Artin braid group
 $\op{Art}_{\Gamma}$ but also a presentation of  $\op{Art}_{\Gamma}$.
 The pair $(\op{Art}_{\Gamma}, \{t_{i} \}_{i \in I})$ consisting of an
 abstract Artin braid group together 
with a set of standard generators is called an {\em Artin
system}. 

\smallskip

\noindent
(ii) Given an Artin system $(\op{Art}_{\Gamma}, \{t_{i}
| i \in I \})$, one obtains a natural character 
\[
\deg : \op{Art}_{\Gamma} \to {\mathbb Z},
\]
which sends each generator to $1 \in {\mathbb Z}$.

\smallskip

\noindent
(iii) If $\Lambda \subset \Gamma$ is a full subgraph (i.e. $I(\Lambda)
\subset I(\Gamma)$ and if $i, j \in I(\Lambda)$, then $k_{ij}(\Lambda)
= k_{ij}(\Gamma)$, then the natural homomorphism $\op{Art}_{\Lambda}
\subset \op{Art}_{\Gamma}$ is known to be injective \cite{vanderlek}.

\end{rem}

\

\bigskip

Let now $(\op{Art}_{\Gamma},\{ t_{i} \}_{i \in I})$ be the Artin
system corresponding to a graph $\Gamma$. We are interested in  the
central extensions of $\op{Art}_{\Gamma}$ by ${\mathbb Z}$ or
equivalently in the group cohomology
$H^{2}(\op{Art}_{\Gamma},{\mathbb Z})$, where ${\mathbb Z}$ is taken
with the trivial $\op{Art}_{\Gamma}$-action. The group
$H^{2}(\op{Art}_{\Gamma},{\mathbb Z})$ can be quite complicated. For
example, due to the work of Arnold \cite{arnold-braid}, Cohen
\cite{cohen-braids} and Fuks
\cite{fuks-braid} it is known that when $\Gamma$ is a Dynkin graph of
type $A_{n}$,  the group $H^{2}(\op{Art}_{\Gamma},{\mathbb
Z})$ is torsion and that when $n \geq 3$ it does contain a 
non-trivial two torsion.

Let ${\mathcal A}(\Gamma)$ be the set of all abelian subgroups
$\op{Art}_{\Gamma}$ which are generated by $\{t_{i}\}_{i \in
J}$ for some $J \subset I$. Equivalently ${\mathcal A}(\Gamma)$ can be
identified with the set of all subgroups $G \subset \op{Art}_{\Gamma}$
of the form $G = \op{Art}_{\Lambda}$, where $\Lambda \subset \Gamma$ is a
full subgraph with no edges. 

Let $\gamma \in H^{2}(\op{Art}_{\Gamma},{\mathbb Z})$ and let 
\begin{equation*}
0 \to {\mathbb Z} \to \Phi_{\gamma} \to \op{Art}_{\Gamma} \to 1,
\label{eq-ext} \tag{$\gamma$}
\end{equation*}
be the corresponding central extension of $\op{Art}_{\Gamma}$. 

Consider the subgroup 
\[
E(\Gamma) := \bigcap_{G \in {\mathcal A}(\Gamma)} \ker\left[ 
H^{2}(\op{Art}_{\Gamma},{\mathbb Z}) \to H^{2}(G,{\mathbb Z})\right] 
\subset H^{2}(\op{Art}_{\Gamma},{\mathbb Z}).
\]
Explicitly $E(\Gamma)$ consists of all $\gamma \in
H^{2}(\op{Art}_{\Gamma},{\mathbb Z})$ for which the natural pullback
sequence 
\[
\xymatrix@R=8pt{
0 \ar[r] & {\mathbb Z} \ar[r] \ar@{=}[d] & \Phi_{\gamma} \ar[r] &
\op{Art}_{\Gamma} \ar[r] & 1 \\
0 \ar[r] & {\mathbb Z} \ar[r] &
\Phi_{\gamma}\times_{\op{Art}_{\Gamma}} G \ar[r] \ar@{^{(}->}[u] & G
\ar[r] \ar@{^{(}->}[u] & 1
}
\]
is split.

We have the following simple

\begin{lem} \label{lem-vanish}
Assume that the graph $\Gamma$ is simply-laced (and hence simply
connected). Then $E(\Gamma) = 0$.
\end{lem}
{\bf Proof.} Fix $\gamma \in E(\Gamma)$. Consider the central
extension \eqref{eq-ext} and let $\{a_{i}\}_{i \in I} \subset
\Phi_{\gamma}$ be lifts of $t_{i} \in \op{Art}_{\Gamma}$.

Since $\Gamma$ is assumed to be simply-laced, it follows that all
relations defining $\op{Art}_{\Gamma}$ are:
\begin{itemize}
\item $t_{i}t_{j} = t_{j}t_{i}$ if $i$ and $j$ are not connected by an
edge;
\item $t_{i}t_{j}t_{i} = t_{j}t_{i}t_{j}$ if $i$ and $j$ are connected
by an edge.
\end{itemize}
Let $i \neq j$ be two vertices of $\Gamma$ which are not connected by
an edge. Consider the subgroup $G = \langle t_{i}, t_{j} \rangle
\subset \op{Art}_{\Gamma}$. Then $G \in {\mathcal A}(\Gamma)$ and so by
our hypothesis this implies that the sequence 
\[
0 \to {\mathbb Z} \to \Phi_{\gamma}\times_{\op{Art}_{\Gamma}} G \to G
\to 1
\]
is split. In particular this means that
$\Phi_{\gamma}\times_{\op{Art}_{\Gamma}} G$ is abelian and so for
$a_{i}, a_{j} \in \Phi_{\gamma}\times_{\op{Art}_{\Gamma}} G \subset
\Phi_{\gamma}$ we get $a_{i}a_{j} = a_{j}a_{i}$.

Let now $c \in \Phi_{\gamma}$ be the generator of  ${\mathbb Z}
\subset \Phi_{\gamma}$ and let $i, j \in I$ be two vertices of $\Gamma$
which are connected by an edge. Since $t_{i}t_{j}t_{i} =
t_{j}t_{i}t_{j}$ we have that
\begin{equation} \label{eq-nij}
a_{i}a_{j}a_{i} = a_{j}a_{i}a_{j}c^{n_{[ij]}},
\end{equation}
for some integer $n_{[ij]}$.

Consider the one dimensional complex $\Gamma$. Let
$C_{1}'(\Gamma,{\mathbb Z})$ be the free abelian group generated by
the oriented edges of $\Gamma$. In particular, for every edge of
$\Gamma$ we have two generators of $C_{1}'(\Gamma,{\mathbb
Z})$. Introduce the relation that the two generators corresponding to
an edge are negative of each other. Let $C_{1}(\Gamma,{\mathbb Z})$
denote the quotient group. It has one generator for each edge of
$\Gamma$. Note that these generators can be denoted by their end
points. We put $[ij]$ for the edge connecting $i$ and $j$, with the
orientation `from $i$ to $j$'. In particular $[ji] = - [ij]$ in
$C_{1}(\Gamma,{\mathbb Z})$. 

The group of 1-cochains of $\Gamma$ with coefficients in  ${\mathbb
Z}$ is the group $\op{Hom}_{{\mathbb Z}}(C_{1}(\Gamma,{\mathbb
Z}),{\mathbb Z})$. In other words, a 1-cochain of $\Gamma$ is given by 
a collection of integers 
\[
\{ f_{[ij]} \in {\mathbb Z} \}_{[ij] 
\text{ is an oriented edge of $\Gamma$ }},
\] 
such that $f_{[ij]} = - f_{[ji]}$.

Note that $n_{[ij]} = - n_{[ji]}$ due to the defining relation
\eqref{eq-nij} and so $n := \{ n_{[ij]} \} \in C^{1}(\Gamma,{\mathbb
Z})$. However $\dim \Gamma =  1$ and so $ C^{1}(\Gamma,{\mathbb
Z}) = Z^{1}(\Gamma,{\mathbb Z})$. Furthermore $\Gamma$ is simply
connected and so $Z^{1}(\Gamma,{\mathbb Z}) = \delta C^{0}(\Gamma,{\mathbb
Z})$. Hence we can find a zero cochain $m := \{ m_{i} \}_{i \in I}$ of
the simplicial complex $\Gamma$, so that $n = \delta m$. 

Consider the elements $b_{i} = a_{i}c^{m_{i}} \in
\Phi_{\Gamma}$. Clearly the $b_{i}$'s also lift the $t_{i}$'s and we
have $b_{i}c = cb_{i}$ and $b_{i}b_{j} = b_{j}b_{i}$ for $i, j \in I$
which are not connected by an edge in $\Gamma$. Finally, for $i, j \in
I$ which are connected by an edge, we calculate
\[
\begin{split}
b_{i}b_{j}b_{i} & = a_{i}a_{j}a_{i}c^{2m_{i} + m_{j}} =
a_{j}a_{i}a_{j}c^{n_{[ij]} + 2m_{i} + m_{j}} = \\
& = b_{j}b_{i}b_{j}c^{n_{[ij]} + 2m_{i} + m_{j} - 2m_{j} + m_{i}} = 
b_{j}b_{i}b_{j}c^{n_{[ij]} + m_{i} - m_{j}} =\\
& = b_{j}b_{i}b_{j}c^{(n - \delta m)_{[ij]}} = b_{j}b_{i}b_{j}.
\end{split}
\]
This implies that the subgroup of $\Phi_{\gamma}$ generated by the
$b_{i}$'s is isomorphic to $\op{Art}_{\Gamma}$ and splits off as a
direct summand in $\Phi_{\gamma}$. Hence $\gamma = 0$ in
$H^{2}(\op{Art}_{\Gamma},{\mathbb Z})$ and so the lemma is
proven. \hfill $\Box$

\section{Displacement angles} \label{sec-angles} 

Let $H$ be a free abelian group of rank $2g$ and let $\theta :
H\otimes H \to {\mathbb Z}$ be a symplectic unimodular pairing on
$H$. Consider the $2g$ dimensional vector space $H_{{\mathbb R}} :=
H\otimes {\mathbb R}$. The real symplectic group $Sp(H_{{\mathbb
R}},\theta)$ is homotopy equivalent to its maximal compact subgroup
which in turn is isomorphic to the unitary group $U(g)$. In particular
$\pi_{1}(Sp(H_{{\mathbb R}},\theta)) \cong \pi_{1}(U(g)) \cong
{\mathbb Z}$ and so the universal cover $\widetilde{Sp}(H_{{\mathbb
R}},\theta)$ of $Sp(H_{{\mathbb R}},\theta)$ is naturally a central
extension
\begin{equation} \label{eq-centralR}
0 \to {\mathbb Z} \to \widetilde{Sp}(H_{{\mathbb
R}},\theta) \to Sp(H_{{\mathbb R}},\theta) \to 1.
\end{equation}
Let $\Lambda(H_{{\mathbb R}},\theta)$ be the Lagrangian Grassmanian of
the symplectic vector space $(H_{{\mathbb R}},\theta)$. The
Grassmanian $\Lambda(H_{{\mathbb R}},\theta)$ can be identified with
the homogeneous space $U(g)/O(g)$ as follows. Choose a complex
structure $I : H_{{\mathbb R}} \to H_{{\mathbb R}}$ which is 
$\theta$-tamed. This simply means that $\gamma(x,y) := \theta(I(x),y)$ is a
positive definite symmetric form and so $\eta = \gamma + \sqrt{-1}\theta$
is a positive definite Hermitian form on the $g$-dimensional complex
vector space $H^{\mathbb C} := (H_{{\mathbb R}},I)$. Now every element  
in the unitary group $U(H^{\mathbb C},\eta)$ necessarily preserves
$\theta$ and so we get an inclusion $U(g) \cong U(H^{\mathbb C},\eta)
\subset  Sp(H_{{\mathbb R}},\theta)$ which is a homotopy equivalence. 
If $\lambda \subset H_{{\mathbb R}}$ is a Lagrangian subspace, then
every basis of $\lambda$ which is orthonormal w.r.t. $\gamma_{|\lambda}$
will also be a ${\mathbb C}$-basis of $H^{\mathbb C}$ which is
orthonormal w.r.t. $\eta$. In particular if $\lambda, \mu \in
\Lambda(H_{{\mathbb R}},\theta)$ are two Lagrangian subspaces and we
choose $\gamma$-orthonormal bases in $\lambda$ and $\mu$ respectively, then
there will be unique element $u \in U(H^{\mathbb C},\eta)$ which sends
the basis for $\lambda$ to the basis for $\mu$ and so $u(\lambda) =
\mu$. This shows that $U(H^{\mathbb C},\eta)$ will act transitively on 
$\Lambda(H_{{\mathbb R}},\theta)$ and that the stabilizer of a point
$\lambda \in \Lambda(H_{{\mathbb R}},\theta)$ in $U(H^{\mathbb
C},\eta)$ can be identified with the orthogonal group $O(\lambda,
\gamma_{|\lambda})$. Thus $\Lambda(H_{{\mathbb R}},\theta) =
U(H^{\mathbb C},\eta)/O(\lambda, \gamma_{|\lambda}) \cong U(g)/O(g)$.

This homogeneous space interpretation can be used to show \cite{arnold-maslov}
that the fundamental group of the Lagrangian Grassmanian is isomorphic to
${\mathbb Z}$. Indeed the natural determinant homomorphism $\det :
U(g) \to S^{1}$ restricts to $\det : O(g) \to \{\pm 1\}$ on $O(g)$ and
so descends to a well defined map $d : U(g)/O(g) \to S^{1}/\{\pm 1\}
\cong S^{1}$. The fiber of $d$ is diffeomorphic to the homogeneous
space $SU(g)/SO(g)$. But $SU(g)$ is simply connected and $SO(g)$ is
connected and so $\pi_{1}(SU(g)/SO(g)) = \{1\}$ from the long exact
sequence of homotopy groups for the fibration $SU(g) \to
SU(g)/SO(g)$. Therefore by the long exact sequence of homotopy groups for
the fibration $d : U(g)/O(g) \to S^{1}$ we conclude that $d$ induces
an isomorphism on fundamental groups,
i.e. $\pi_{1}(\Lambda(H_{{\mathbb R}},\theta)) = \pi_{1}(U(g)/O(g)) =
\{1\}$. 

Note that the group
$Sp(H_{{\mathbb R}},\theta)$ also acts transitively on $\Lambda(H_{{\mathbb
R}},\theta)$ and that $\widetilde{Sp}(H_{{\mathbb R}},\theta)$ acts
transitively on the universal cover $\widetilde{\Lambda}(H_{{\mathbb
R}},\theta)$ of $\Lambda(H_{{\mathbb R}},\theta)$. 

Recall that every vector $a \in H_{{\mathbb R}}$ generates a one
parameter unipotent subgroup in $Sp(H_{{\mathbb
R}},\theta)$ by the formula
\[
\xymatrix@R=4pt{
T_{a} : & {\mathbb R} \ar[r] & Sp(H_{{\mathbb R}},\theta) \\
& s \ar[r] & {\protect (x \mapsto x + s\theta(a,x)a)}.
}
\]
Every element of the form $T_{a}(s) = T_{\sqrt{s}a}(1)$ is called a {\em
symplectic transvection}. In the case when $H = H_{1}(C,{\mathbb Z})$
for some smooth surface $C$, the element $T_{a}(1)$ is the image of
the oriented Dehn twist along a simple closed curve representing the
homology class $a \in H_{1}(C,{\mathbb Z})$.

We begin with the following lemma.

\begin{lem} \label{lem-lift} Let $t \in Sp(H_{{\mathbb R}},\theta)$ be
a symplectic transvection. Then there exists a unique lift $\tilde{t}
\in \widetilde{Sp}(H_{{\mathbb R}},\theta)$ of $t$ which acts with
fixed points on $\widetilde{\Lambda}(H_{{\mathbb R}},\theta)$.
\end{lem}
{\bf Proof.} To check that $\tilde{t}$ exist write $t = T_{a}(1)$ for
some $a \in H_{{\mathbb R}}$. The vector $a$ can be included in a
symplectic basis $a_{1} = a, a_{2}, \ldots, a_{g}, b_{1}, b_{2},
\ldots, b_{g}$ of $H_{{\mathbb R}}$, and so $t$ preserves the Lagrangian
subspace $\lambda := \op{Span}_{{\mathbb R}}(a_{1},\ldots,
a_{g})$. Let $ \tilde{\lambda} \in
\widetilde{\Lambda}(H_{{\mathbb R}},\theta)$ be a preimage of $\lambda \in
\Lambda(H_{{\mathbb R}},\theta)$ and let $p \in
\widetilde{Sp}(H_{{\mathbb R}},\theta)$ be a preimage of $t \in
Sp(H_{{\mathbb R}},\theta)$. Then
$p\cdot \tilde{\lambda}$ maps to $t\cdot \lambda = \lambda$ and so we can find
a  deck transformation $c \in 
\pi_{1}({\Lambda}(H_{{\mathbb R}},\theta),\lambda) = {\mathbb Z}$
satisfying $c(p \cdot \tilde{\lambda}) = \tilde{\lambda}$. On the other
hand, the continuous map $m : Sp(H_{{\mathbb R}},\theta) \to
{\Lambda}(H_{{\mathbb R}},\theta)$, given by $m(g) = g\cdot \lambda$
induces a homomorphism $m_{*} : \pi_{1}(Sp(H_{{\mathbb R}},\theta),e)
\to \pi_{1}({\Lambda}(H_{{\mathbb R}},\theta),\lambda)$. If $c = m_{*}(\tilde{c})$ for some $\tilde{c} \in
\pi_{1}(Sp(H_{{\mathbb R}},\theta),e) \subset
Z(\widetilde{Sp}(H_{{\mathbb R}},\theta))$, we get
\[
\tilde{\lambda} = c(p\cdot \tilde{\lambda}) =
m_{*}(\tilde{c})(p\cdot \tilde{\lambda}) = (\tilde{c}p)\cdot \tilde{\lambda},
\]
where by $(\tilde{c}p)$ we mean the product of $\tilde{c} \in
Z(\widetilde{Sp}(H_{{\mathbb R}},\theta)) \subset
\widetilde{Sp}(H_{{\mathbb R}},\theta)$ and $p \in
\widetilde{Sp}(H_{{\mathbb R}},\theta)$ in the group
$\widetilde{Sp}(H_{{\mathbb R}},\theta)$. However multiplication by
elements in $\pi_{1}(Sp(H_{{\mathbb R}},\theta),e) \subset
Z(\widetilde{Sp}(H_{{\mathbb R}},\theta))$ preserves the fibers of the
covering map $\widetilde{Sp}(H_{{\mathbb R}},\theta) \to
Sp(H_{{\mathbb R}},\theta)$. and so we may take $\tilde{t} :=
\tilde{c}p$.

Therefore, in order to finish the proof of the existence of
$\tilde{t}$ we have to show that $c \in \op{im}(m_{*})$. As
we explained above the identification $\Lambda(H_{{\mathbb R}},\theta)
= U(H^{\mathbb C},\eta)/O(\lambda,\gamma_|\lambda)$ implies that the
map $m_{*}$ fits in the following commutative diagram with exact rows
\[
\xymatrix{
0 \ar[r] & \pi_{1}(Sp(H_{{\mathbb R}},\theta),e) \ar[r]^{m_{*}}
\ar@{=}[d] & \pi_{1}({\Lambda}(H_{{\mathbb R}},\theta),\lambda) \ar[r]
\ar@{=}[d] & {\mathbb Z}/2 \ar[r] \ar@{=}[d] & 0 \\
0 \ar[r] & {\mathbb Z} \ar[r]_{\op{mult}_{2}} & {\mathbb Z} \ar[r] & 
{\mathbb Z}/2 \ar[r] & 0.}
\]
In particular if $\Lambda_{(2)}(H_{{\mathbb R}},\theta) \to 
\Lambda(H_{{\mathbb R}},\theta)$  denotes
 the unramified double cover corresponding to the
surjection $\pi_{1}({\Lambda}(H_{{\mathbb R}},\theta),\lambda) \to
{\mathbb Z}/2$, then $m$ factors as
$Sp(H_{{\mathbb R}},\theta) \to \Lambda_{(2)}(H_{{\mathbb R}},\theta) \to 
\Lambda(H_{{\mathbb R}},\theta)$. Furthermore, by the definition of
$m$ we have $m(tg) = t\cdot m(g)$ for all $g \in Sp(H_{{\mathbb
R}},\theta)$ and so we have a natural lift $t_{(2)}$ of the action of $t$ on
$\Lambda(H_{{\mathbb R}},\theta)$ to an automorphism of the double
cover $\Lambda_{(2)}(H_{{\mathbb R}},\theta)$. Since the map 
$Sp(H_{{\mathbb R}},\theta) \to \Lambda_{(2)}(H_{{\mathbb R}},\theta)$
has connected fibers, it induces an isomorphism on fundamental groups
and so $c \in \op{im}(m_{*})$ if and only if the automorphism
$t_{(2)}$ acts trivially on the fiber of $\Lambda_{(2)}(H_{{\mathbb
R}},\theta) \to  \Lambda(H_{{\mathbb R}},\theta)$ over the point
$\lambda \in \Lambda(H_{{\mathbb R}},\theta)$. But by construction 
this fiber can be identified with the quotient ${\mathbb Z}/2$ of
$\pi_{1}(\Lambda(H_{{\mathbb R}},\theta),\lambda)$ and the action 
$t_{(2)}$ on the fiber can be identified with the action on ${\mathbb
Z}/2$ induced from $t_{*} : \pi_{1}(\Lambda(H_{{\mathbb
R}},\theta),\lambda) \to \pi_{1}(\Lambda(H_{{\mathbb
R}},\theta),\lambda)$. Finally note that $Sp(H_{{\mathbb R}},\theta)$
is connected and so the action of $t$ on  $\Lambda(H_{{\mathbb
R}},\theta)$ is homotopic to the identity. Hence $t_{*}$ acts
trivially on $\pi_{1}(\Lambda(H_{{\mathbb R}},\theta),\lambda)$ and
$c \in \op{im}(m_{*})$.

Next we prove the uniqueness of $\tilde{t}$. To that end we calculate
the fixed locus of $\tilde{t}$ on $\widetilde{\Lambda}(H_{{\mathbb
R}},\theta)$ explicitly. Clearly the image of
$\op{Fix}_{\tilde{t}}(\widetilde{\Lambda}(H_{{\mathbb R}},\theta))$ in
$\Lambda(H_{{\mathbb R}},\theta)$ is contained in the fixed locus 
$\op{Fix}_{t}(\Lambda(H_{{\mathbb R}},\theta))$. Since
$\pi_{1}(Sp(H_{{\mathbb R}},\theta))$ is central in
$\widetilde{Sp}(H_{{\mathbb R}},\theta)$, it follows that the element
$\tilde{t} \in \widetilde{Sp}(H_{{\mathbb R}},\theta)$ commutes with
all $c \in \pi_{1}(Sp(H_{{\mathbb R}},\theta))$. Since
$\pi_{1}(Sp(H_{{\mathbb R}},\theta))$ acts transitively on the fibers
of $\pi : \widetilde{\Lambda}(H_{{\mathbb R}},\theta) \to
\Lambda(H_{{\mathbb R}},\theta)$ and $\tilde{t}\cdot \tilde{\lambda} =
\tilde{\lambda}$, it  follows that $\tilde{t}$ fixes all points in 
$\widetilde{\Lambda}(H_{{\mathbb R}},\theta)$ that map to $\lambda \in 
\Lambda(H_{{\mathbb R}},\theta)$. In particular $\pi^{-1}(\lambda)
\subset \op{Fix}_{\tilde{t}}(\widetilde{\Lambda}(H_{{\mathbb
R}},\theta))$. In fact we have
\begin{equation} \label{eq-tildet-fixed}
\op{Fix}_{\tilde{t}}(\widetilde{\Lambda}(H_{{\mathbb R}},\theta)) =
\pi^{-1}(\op{Fix}_{t}(\Lambda(H_{{\mathbb R}},\theta))).
\end{equation}
Indeed, the locus $\op{Fix}_{t}(\Lambda(H_{{\mathbb R}},\theta))$
consists of all Lagrangian subspaces $\mu \in 
\Lambda(H_{{\mathbb R}},\theta)$ that contain the vector $a \in
H_{{\mathbb R}}$ and so is isomorphic to the Lagrangian Grassmanian of
a symplectic vector space of dimension $2g - 2$. To see this consider
the $\theta$-orthogonal complement $a^{\perp}$ of $a$ in $H_{{\mathbb
R}}$. The $2g-1$ dimensional subspace $a^{\perp} \subset H_{\mathbb R}$
contains $a$ and inherits a skew-symmetric form $\theta_{|a^{\perp}}$
whose kernel is spanned by $a$. The $\theta$-lagrangian subspaces in  
$H_{{\mathbb R}}$ which contain the vector $a$ are all contained in
$a^{\perp}$. Therefore $\op{Fix}_{t}(\Lambda(H_{{\mathbb R}},\theta))$
can be identified with the set of all $g$-dimensional subspaces in
$a^{\perp}$ which are $\theta$-isotropic and contain the vector
$a$. Since $\ker(\theta_{|a^{\perp}}) = {\mathbb R}\cdot a$, the form 
$\theta_{|a^{\perp}}$ descends to a symplectic form $\bar{\theta}$ on
the quotient space $\overline{H}_{{\mathbb  R}} := a^{\perp}/{\mathbb
R}\cdot a$.  Now every $g$-dimensional $\theta$-isotropic subspace in
$a^{\perp}$ which contains $a$ will map onto a
$\bar{\theta}$-Lagrangian subspace of $\overline{H}_{{\mathbb R}}$ and
conversely - the preimage of a $\bar{\theta}$-Lagrangian subspace of
$\bar{\mu} \subset \overline{H}_{{\mathbb R}}$ will be a
$g$-dimensional $\theta$-isotropic subspace $\mu \subset a^{\perp}$
which contains $a$. In other words, we have constructed an inclusion 
$\Lambda(\overline{H}_{{\mathbb R}},\bar{\theta}) \hookrightarrow
\Lambda(H_{{\mathbb R}},\theta)$, whose image is precisely
$\op{Fix}_{t}(\Lambda(H_{{\mathbb R}},\theta))$. In particular this
shows that $\op{Fix}_{t}(\Lambda(H_{{\mathbb R}},\theta))$ is
connected and that the inclusion map $\Lambda(\overline{H}_{{\mathbb
R}},\bar{\theta}) = \op{Fix}_{t}(\Lambda(H_{{\mathbb R}},\theta))
\hookrightarrow \Lambda(H_{{\mathbb R}},\theta)$ induces an
isomorphism on fundamental groups. This implies that
$\pi^{-1}(\op{Fix}_{t}(\Lambda(H_{{\mathbb R}},\theta)))$ is
and simply connected and that the map $\pi :
\pi^{-1}(\op{Fix}_{t}(\Lambda(H_{{\mathbb R}},\theta))) \to
\op{Fix}_{t}(\Lambda(H_{{\mathbb R}},\theta))$ is the universal
covering map for $\op{Fix}_{t}(\Lambda(H_{{\mathbb R}},\theta))$. But
by definition $\pi\circ \tilde{t} = t\circ\pi$ and that $t$ acts
trivially on $\op{Fix}_{t}(\Lambda(H_{{\mathbb R}},\theta))$. This
shows that
$\tilde{t}$ is an automorphism of the universal covering map 
$\pi : \pi^{-1}(\op{Fix}_{t}(\Lambda(H_{{\mathbb R}},\theta))) \to
\op{Fix}_{t}(\Lambda(H_{{\mathbb R}},\theta))$ and so
$\tilde{t|\pi^{-1}(\op{Fix}_{t}(\Lambda(H_{{\mathbb R}},\theta)))}$
must be a deck transformation. However by construction $\tilde{t}$
fixes $\tilde{\lambda} \in \pi^{-1}(\op{Fix}_{t}(\Lambda(H_{{\mathbb
R}},\theta)))$ and so
$\tilde{t|\pi^{-1}(\op{Fix}_{t}(\Lambda(H_{{\mathbb R}},\theta)))} =
\op{id}$. Thus we have established the validity of
\eqref{eq-tildet-fixed}.

Consider now some other lift $q$ of $t$. Then $q = c\tilde{t}$ for
some $c \in \pi_{1}(Sp(H_{{\mathbb R}},\theta))$. Now again 
$\op{Fix}_{q}(\widetilde{\Lambda}(H_{{\mathbb R}},\theta)) \subset 
\pi^{-1}(\op{Fix}_{t}(\Lambda(H_{{\mathbb R}},\theta))$ and since
$\tilde{t}$ acts trivially on
$\pi^{-1}(\op{Fix}_{t}(\Lambda(H_{{\mathbb R}},\theta))$ we see that
$q_{|\pi^{-1}(\op{Fix}_{t}(\Lambda(H_{{\mathbb R}},\theta))} =
c_{|\pi^{-1}(\op{Fix}_{t}(\Lambda(H_{{\mathbb R}},\theta))}$. But $c$
is a deck transformation and so acts without fixed points. Thus
$\op{Fix}_{q}(\widetilde{\Lambda}(H_{{\mathbb R}},\theta)) =
\varnothing$ and the uniqueness of the lift $\tilde{t}$ is
proven. \hfill $\Box$

\bigskip

\begin{figure}[!ht]
\begin{center}
\epsfig{file=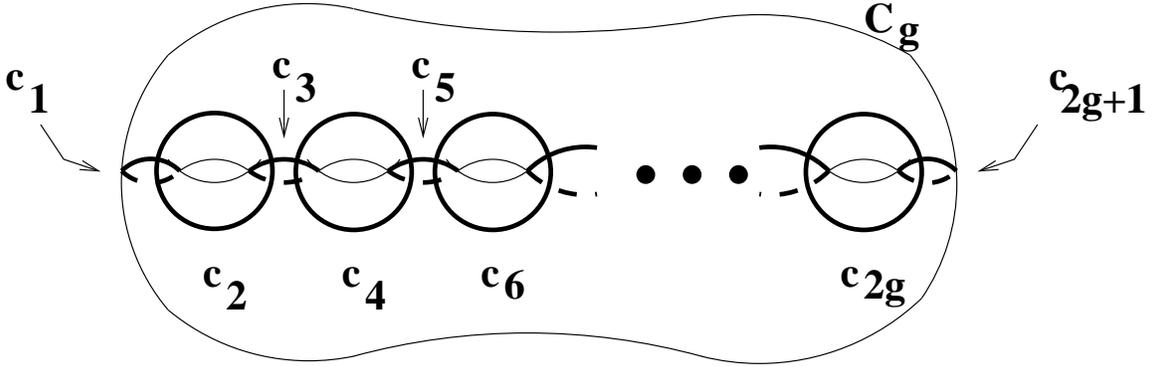,width=6in} 
\end{center}
\caption{Loops generating $B_{2g+2}$.}\label{fig-generators} 
\end{figure}

Consider next the standard Artin braid group $B_{2g + 2}$ on $2g+2$
strands. In other words $B_{2g + 2} :=
\op{Art}_{\boldsymbol{A}_{2g+1}}$ is the Artin braid group
corresponding to the Dynkin graph  $\boldsymbol{A}_{2g+1}$. Explicitly
\[
B_{2g + 2} = \left\langle t_{1}, t_{2}, \ldots, t_{2g+1} \left|
t_{i}t_{j} = t_{j}t_{i} \text{ for } |i - j| \geq 2, t_{i}t_{i+1}t_{i}
= t_{i+1}t_{i}t_{i+1}
\right.\right\rangle.
\]
Fix a closed oriented surface $C_{g}$ of genus $g$ and a 
hyperelliptic involution $\iota$ on $C_{g}$. If we
choose a sequence of loops $c_{1}, c_{2}, \ldots, c_{2g+1}$ on $C_{g}$
as depicted on Figure~\ref{fig-generators}, we can realize the
generators $t_{i}$ geometrically as the right handed Dehn twists
$t_{c_{i}}$. The assignment $t_{i} \to t_{c_{i}} \in \op{Map}_{g,1}$ induces
a homomorphism $\kappa_{g,1} : B_{2g+2} \to \op{Map}_{g,1}$ and after
compositions with the natural projections induces homomorphisms 

\[
\xymatrix@R=4pt{
\kappa_{g}^{1} : & B_{2g+2} \ar[r] &  \op{Map}_{g}^{1}, \\
\kappa_{g} : & B_{2g+2} \ar[r] & \op{Map}_{g}, \\ 
\sigma_{g} : & B_{2g+2} \ar[r] & 
\op{Sp}(H_{1}(C_{g},{\mathbb Z})).
}
\]
Furthermore, if one is
careful enough to chose the $c_{i}$ so that they are invariant under the
hyperelliptic involution $\iota$, then the image of $\kappa_{g,1}$
will be the hyperelliptic mapping class group $\Delta_{g,1}$.
It is known \cite{birman-hilden} that $\kappa_{g,1} : B_{2g+2} \to
\Delta_{g,1}$ is surjective with a kernel normally generated by the
element $(t_{1}\ldots t_{2g+1})^{2g+1}(t_{1}t_{2}\ldots
t_{2g+1}t_{2g+1} \ldots t_{2}t_{1})^{-1}$, and that $\kappa_{g}^{1} :
B_{2g+2} \to \Delta_{g}^{1}$ is surjective with a kernel normally
generated by the elements $(t_{1}\ldots t_{2g+1})^{2g+1}(t_{1}t_{2}\ldots
t_{2g+1}t_{2g+1} \ldots t_{2}t_{1})^{-1}$ and $(t_{1}\ldots
t_{2g+1})^{2g+2}$. It is also known \cite{varchenko-braid} that the
map $\sigma_{g}$ is \label{kappa}
surjective, and so we have a sequence of surjective group
homomorphisms
\[
B_{2g+2} \twoheadrightarrow \Delta_{g,1} \twoheadrightarrow
\Delta_{g}^{1} \twoheadrightarrow \Delta_{g} \twoheadrightarrow
\op{Sp}(H_{1}(C_{g},{\mathbb Z})).
\]
Consider the lattice $H :=  H_{1}(C_{g},{\mathbb Z})$ together with
symplectic unimodular pairing $\theta : H\otimes H \to {\mathbb Z}$
corresponding to the intersection of cycles. Let now 
\begin{equation} \label{eq-braid-extension}
0 \to {\mathbb Z} \to \widetilde{B}_{2g+2} \to B_{2g+2} \to 0
\end{equation}
be the pullback of the central extension \eqref{eq-centralR} via the
homomorphism 
\[
\sigma_{g} : B_{2g +2} \to Sp(H,\theta) \subset
Sp(H_{{\mathbb R}},\theta).
\]
we have the following important
\begin{prop} \label{prop-split}
The extension \eqref{eq-braid-extension} is a split extension.
\end{prop} 
{\bf Proof.} Let $\gamma \in H^{2}(B_{2g+2},{\mathbb Z}) =
H^{2}(\op{Art}_{\boldsymbol{A}_{2g+2}},{\mathbb Z})$ be the
extension class of \eqref{eq-braid-extension}. In view of
Lemma~\ref{lem-vanish} it suffices to show that $\gamma \in
E(\boldsymbol{A}_{2g+2})$ or equivalently that $\gamma$ splits when
restricted on every $G \in {\mathcal A}(\boldsymbol{A}_{2g+2})$. Since
by definition \eqref{eq-braid-extension} is a pullback of a central
extension of the group $Sp(H,\theta)$, it suffices to check that the
extension \eqref{eq-centralR} splits when restricted on any abelian subgroup
of $Sp(H,\theta)$ which is generated by the Dehn twists about any
finite collection of non-intersecting $c_{i}$'s. 

Let $\{a_{1}, \ldots, a_{k} \} \subset \{c_{1}, \ldots, c_{2g+1}\}$ be
such that $\theta(a_{i},a_{j}) = 0$ for all $i, j$. Let $t_{i} :=
T_{a_{i}}(1) \in Sp(H,\theta)$ be the corresponding symplectic
transvections and let $S \subset Sp(H,\theta)$ be the subgroup
generated by the $t_{i}$'s. Consider the pullback of the extension 
\eqref{eq-centralR} via the inclusion map $S \hookrightarrow
Sp(H,\theta) \subset Sp(H_{{\mathbb R}},\theta)$:
\begin{equation} \label{eq-restrict}
0 \to {\mathbb Z} \to \widetilde{S} \to S \to 0.
\end{equation}
As explained above, it suffices to show that \ref{eq-restrict} is split 
in order to prove the proposition. To achieve this consider the
subspace $\op{Span}_{{\mathbb R}}(a_{1}, \ldots, a_{k}) \subset
H_{{\mathbb R}}$. Since by assumption this subspace is
$\theta$-isotropic we can find a Lagrangian subspace $\lambda \in
\Lambda(H_{{\mathbb R}},\theta)$ so that $\lambda \supset 
\op{Span}_{{\mathbb R}}(a_{1}, \ldots, a_{k})$. In particular we will
have $t_{i}\cdot \lambda = \lambda$ for all $i = 1, \ldots, k$ and
hence $\lambda \in \op{Fix}_{S}(\Lambda(H_{{\mathbb R}},\theta)) \neq
\varnothing$. We now have the following

\begin{lem} \label{lem-lift-S} For every element $g \in S$ there
exists a lift $\tilde{g} \in \widetilde{S}$ which is uniquely
cha\-ra\-cte\-ri\-zed by the property
$\op{Fix}_{\tilde{g}}(\widetilde{\Lambda}(H_{{\mathbb R}},\theta)
\supset \pi^{-1}(\op{Fix}_{S}(\Lambda(H_{{\mathbb R}},\theta)))$. 
\end{lem}
{\bf Proof.} We will prove the lemma by induction on $k$.  
For $k = 1$ this is precisely the statement of Lemma~\ref{lem-lift}. 
Let $k > 1$. By the inductive 
hypothesis the elements of the 
subgroup generated by $t_{1}, \ldots, t_{k-1}$ lift uniquely to
elements in $\widetilde{S}$ which fix all points in the set 
$\pi^{-1}(\op{Fix}_{\langle t_{1}, \ldots, t_{k-1}
\rangle}(\Lambda(H_{{\mathbb R}},\theta)))$. By 
applying Lemma~\ref{lem-lift} again to the image of the 
element $a_{k}$ in the symplectic vector
space $\op{Span}(a_{1}, \ldots, a_{k-1})^{\perp}/\op{Span}(a_{1},
\ldots, a_{k-1})$ we get the required lifts for all elements in
$S$. \hfill $\Box$

\medskip

By Lemma~\ref{lem-lift-S} we get a set theoretic map $S \to
\widetilde{S}$, $g \mapsto \tilde{g}$, which splits the exact sequence 
\eqref{eq-restrict}. However if $g, h \in S$ are two elements, then
$\tilde{g}\cdot \tilde{h}$ is a lift of $g\cdot h$ which necessarily
fixes all points in $\pi^{-1}(\op{Fix}_{S}(\Lambda(H_{{\mathbb
R}},\theta))$ since $\tilde{g}$ and $\tilde{h}$ fix those points
individually. Thus $\widetilde{g\cdot h} = \tilde{g}\cdot \tilde{h}$
and hence the assignment $g \mapsto \tilde{g}$ is a group-theoretic
splitting of \eqref{eq-restrict}. This finishes the proof of
proposition \ref{prop-split}. \hfill $\Box$

\bigskip

The reasoning in the proof of the 
previous proposition gives as an immediate corollary the following
statement which we record here for future use. 
\begin{corr} \label{cor-homomorphism}
There exists a homomorphism  $\tilde{\sigma}_{g} : B_{2g+2} \to
\widetilde{Sp}(H,\theta)$ so that $\pi\circ \sigma_{g} =
\tilde{\sigma}_{g}$ and for which $\tilde{\sigma}_{g}(t_{i}) =
\widetilde{\sigma_{g}(t_{i})}$ is the unique lift of $\sigma_{g}(t_{i})$ 
from Lemma~\ref{lem-lift}.
\end{corr}
\

For the remainder of the paper we fix once and for all a base point
$\lambda_{0} \in \Lambda(H_{{\mathbb 
R}},\theta)$. Let  $U(g) \cong K \subset
Sp(H_{1}(C_{g},{\mathbb R}))$ be a maximal compact subgroup. The
choice of $\lambda_{0}$ determines a $K$-equivariant surjection 
$K \to \Lambda(H_{{\mathbb R}},\theta)$, which
as explained at the beginning of this
section combines with the determinant homomorphism $\det : U(g) \to
S^{1}$ into a well
defined map $d_{K} : \Lambda(H_{{\mathbb R}},\theta) \to
S^{1}$. Let
$\tilde{d}_{K} : \widetilde{\Lambda}(H_{{\mathbb
R}},\theta) \to {\mathbb 
R}$ be a lift of $d_{K}$ to the universal cover. 

Consider the subset 
\[
^{K}\!\widetilde{Sp}(-) := \{ s \in
\widetilde{Sp}(H_{{\mathbb R}},\theta) | 
\tilde{d}_{K}(s\cdot \tilde{\lambda}) \leq
\tilde{d}_{K}(\tilde{\lambda}) 
\text{ for all } \tilde{\lambda} \in \widetilde{\Lambda}(H_{{\mathbb
R}},\theta) \},
\]
and let $\widetilde{Sp}(-) := \cap_{K} \; ^{K}\!\widetilde{Sp}(-)$. 
Clearly $^{K}\!\widetilde{Sp}(-)$ and $\widetilde{Sp}(-)$ are sub
semi-groups in 
$\widetilde{Sp}(H_{{\mathbb R}},\theta)$ and for every simple right
Dehn twist $t \in Sp( H_{{\mathbb R}},\theta)$ we have that $\tilde{t} \in
\widetilde{Sp}(-)$ for the lift $\tilde{t}$ from Lemma~\ref{lem-lift}.

For future reference we give an alternative description of 
the elements in the subsemigroup $\widetilde{Sp}(-)$ in terms of
linear algebraic data. First we have the following easy lemma.

\begin{lem} \label{lem-independence} Let $\widetilde{K}$ be the
universal cover of the maximal compact subgroup $K$. Let
$\widetilde{\det} : \widetilde{K} \to {\mathbb R}$ 
denote the lift of $\det : K\cong U(g) \to S^{1}$ as a group
homomorphism. Then an element
${\mathfrak u} \in \widetilde{K}$ will belong to the subsemigroup
$^{K}\!\widetilde{Sp}(-)$ if and only if $\widetilde{\det}({\mathfrak u})
\leq 0$. In other words
\[
K \cap \, ^{K}\!\widetilde{Sp}(-) = \widetilde{\det}^{-1}({\mathbb
R}_{\leq 0}). 
\]
\end{lem}
{\bf Proof.}  Let $\tilde{\lambda} \in \widetilde{\Lambda}(H_{{\mathbb
R}},\theta)$ and let ${\mathfrak v} \in
\widetilde{K}$ be an element 
which maps to $\tilde{\lambda}$ under the natural map
$\widetilde{K} \to \widetilde{\Lambda}(H_{{\mathbb
R}},\theta)$. Then by the definition of the lifts $\tilde{d}_{K}$ and
$\widetilde{\det}$ we have $\widetilde{\det}({\mathfrak v}) =
\tilde{d}_{K}(\tilde{\lambda}) + c$ where $c \in {\mathbb R}$ is a fixed
constant depending on the choice of the lift $\tilde{d}$ only, and not
on the choices of $\tilde{\lambda}$ or ${\mathfrak v}$. Similarly 
$\widetilde{\det}({\mathfrak u}\cdot
{\mathfrak v}) = \tilde{d}_{K}({\mathfrak u}\cdot \tilde{\lambda}) + c$ for
all ${\mathfrak u} \in \widetilde{K}$. Therefore
\[
\begin{split}
\tilde{d}_{K}(\tilde{\lambda}) -
\tilde{d}_{K}({\mathfrak  u}\cdot\tilde{\lambda}) &  = 
(\widetilde{\det}({\mathfrak v}) - c) - (\widetilde{\det}({\mathfrak u}\cdot
{\mathfrak v}) - c)  \\
& =   \widetilde{\det}({\mathfrak v}) - (\widetilde{\det}({\mathfrak
u}) + \widetilde{\det}({\mathfrak 
v})) = - \widetilde{\det}({\mathfrak
u}).
\end{split}
\]
The lemma is proven. \hfill $\Box$

\medskip

The previous lemma gives a description of the semigroups 
 $K\cap \, ^{K}\!\widetilde{Sp}(-)$ and
$^{K}\!\widetilde{Sp}(-)$ in terms of the natural maps
$\widetilde{\det}$ and $\tilde{d}_{K}$. We would like to have a
similar intrinsic description for the smaller semigroups 
$K\cap \widetilde{Sp}(-)$ and
$\widetilde{Sp}(-)$ respectively.

To that end consider the Lie algebra $\liesp(H_{{\mathbb
R}},\theta) \subset \op{End}(H_{{\mathbb R}})$. Every element $x \in
\liesp(H_{{\mathbb R}},\theta)$ defines a real valued symmetric bilinear form
$\gamma_{x}$ on $H_{{\mathbb R}}$ via the formula
$\gamma_{x}(\bullet,\bullet) := \theta(x(\bullet),\bullet)$. Using the
assignment $x \mapsto \gamma_{x}$ we can now define the {\em
non-positive subcone}  in the
Lie algebra $\liesp(H_{{\mathbb R}},\theta)$ as the cone
\[
\liesp(-) := \{ x \in \liesp(H_{{\mathbb R}},\theta) | \gamma_{x}
\text{ is non-positive definite} \}.
\]
Similarly, for the Lie algebra ${\mathfrak k} \subset
\liesp(H_{{\mathbb R}},\theta)$ of the compact group $K$ we have a
non-positive cone ${\mathfrak k}(-) := {\mathfrak k}\cap \liesp(-)$. 
The elements of this subcone admit a particularly simple
characterization in the fundamental representation of ${\mathfrak k}
\cong {\mathfrak u}(g)$:

\begin{lem} \label{lem-cone} {\em (i)} An element $x \in {\mathfrak k}
\subset \liesp(H_{{\mathbb R}},\theta)$ belongs to the non-positive
subcone ${\mathfrak k}(-)$ if and only if all the eigenvalues of $x$
in the fundamental representation of ${\mathfrak k}
\cong {\mathfrak u}(g)$ are contained in $\sqrt{-1}\cdot {\mathbb R}_{\leq 0}$.

\smallskip

\noindent
{\em (ii)} The interior of the non-positive cone $\liesp(-)$ is the
union of the interiors of all cones of the form ${\mathfrak k}(-)$.
\end{lem}
{\bf Proof.} To prove part (i) recall that the choice of a maximal
compact subgroup $K \subset Sp(H_{{\mathbb R}},\theta)$ corresponds to
the choice of a complex structure $I : H_{{\mathbb R}} \to H_{{\mathbb
R}}$ which is $\theta$-tamed. Once such an $I$ is chosen we can
identify $K$ with the unitary group $U(H^{{\mathbb C}},\eta)$, where 
$H^{{\mathbb C}} = (H_{{\mathbb R}}, I)$ and $\eta(x,y) =
\theta(I(x),y) + \sqrt{-1}\cdot \theta(x,y)$. Under this
identification $H^{{\mathbb C}}$ becomes the fundamental
representation of ${\mathfrak k}$. Furthermore if $a_{1}, \ldots,
a_{g}$ is a basis of $\lambda_{0} \subset H_{{\mathbb R}}$, which is
orthonormal w.r.t. to the form $\theta(I(x),y)$, it follows that
$a_{1}, \ldots, a_{g}$ is an orthonormal basis of $(H^{{\mathbb
C}},\eta)$ and 
$a_{1}, \ldots, a_{g}, I(a_{1}), \ldots, I(a_{g})$ is a symplectic
basis of $(H_{{\mathbb R}},\theta)$. In particular in the basis
$a_{1}, \ldots, a_{g}, I(a_{1}), \ldots, I(a_{g})$ the Gram matrix of
the skew-symmetric form $\theta$ equals the matrix of the linear
operator $I$ equals the standard matrix
\[
\begin{pmatrix}
0 & {\mathbb I}_{g} \\
- {\mathbb I}_{g} & 0
\end{pmatrix}.
\]
(Here as usual ${\mathbb I}_{g}$ denotes the identity $g\times g$ matrix.)
Let now $x \in {\mathfrak k}$ be any element. If we view $x$ as a
linear operator on the complex vector space $H^{{\mathbb C}}$, then in
the basis $a_{1}, \ldots, a_{g}$ this linear operator is given by some
skew-hermitian $g\times g$ matrix $A$. Since every skew-hermitian
matrix can be diagonalized in an orthonormal basis we may assume
without a loss of generality that the basis $a_{1}, \ldots, a_{g}$ is
chosen so that $A = \sqrt{-1}\cdot D$ where $D$ is a  real diagonal
$g\times g$ matrix. Therefore the Gram matrix of the symmetric form
$\gamma_{x}$ in the basis $a_{1}, \ldots, a_{g}, I(a_{1}), \ldots,
I(a_{g})$ is the matrix 
\[
\begin{pmatrix}
0 & D \\
-D & 0
\end{pmatrix}^{t}\cdot
\begin{pmatrix}
0 & {\mathbb I}_{g} \\
- {\mathbb I}_{g} & 0
\end{pmatrix} = 
\begin{pmatrix}
D & 0 \\
0 & D
\end{pmatrix}.
\]
This proves part (i) of the lemma. Part (ii) is straightforward and is
left to the reader. \hfill $\Box$

\bigskip

\begin{rem} \label{rem-nilpotent} \quad {\bf (i)} Since every
unipotent element in the symplectic group $Sp(H_{{\mathbb R}},\theta)$
can be written as a limit of conjugates of elements in a fixed maximal
compact subgroup, part (ii) of the previous lemma implies that the
boundary $\partial \liesp(-)$ of the non-positive cone contains the
intersection ${\mathcal N}\cap \liesp(-)$ of the nilpotent cone 
${\mathcal N} \subset \liesp(H_{{\mathbb R}},\theta)$  and the
non-positive cone. 

\smallskip

\noindent
{\bf (ii)} The same reasoning as in the proof of part (i) of the
previous lemma shows that $\partial \liesp(-)$ consists of all
elements $x \in \liesp(-)$ for which we can find a decomposition 
$H_{{\mathbb R}} = W'\oplus W''$ so that:
\begin{itemize}
\item $W', W'' \subset H_{{\mathbb R}}$ are non-trivial symplectic
subspaces;
\item $x(W') \subset W'$ and $x(W'') \subset W''$;
\item $x_{|W'} \in {\mathcal N}(\liesp(W'))\cap \liesp(W')(-)$ and
$x_{|W''} \in {\mathfrak k}''(-)$ for some maximal compact subgroup
$K'' \subset Sp(W'')$.
\end{itemize}
\end{rem}

\

\noindent
To finish our linear-algebraic description of the semigroups 
$K\cap \widetilde{Sp}(-)$ and $\widetilde{Sp}(-)$ it remains only to
 observe that by Lemmas \ref{lem-independence} and 
\ref{lem-cone} the exponential map $\op{Exp} :
\liesp(H_{{\mathbb R}},\theta) \to \widetilde{Sp}(H_{{\mathbb
R}},\theta)$ maps the cone ${\mathfrak k}(-)$ to the semi-group 
$K\cap \widetilde{Sp}(-)$. Since the exponential map for $U(g)$ is
surjective this implies that the map $\op{Exp} : \liesp(-) \to
\widetilde{Sp}(-)$ is surjective and so $\widetilde{Sp}(-)$ should be
simply thought as the exponentiation of $\liesp(-)$. Therefore the
restriction of $\op{Exp}$ on  the domain 
\[
\liesp(-2\pi,0] := \{ x \in \liesp(H_{{\mathbb R}},\theta) |
\op{spectrum}(\gamma_{x}) \subset (-2\pi,0] \} \subset \liesp(-),
\]
induces a homeomorphism between $\liesp(-2\pi,0]$ and $\widetilde{Sp}$.

Keeping this in mind we can make the following general
definition: 

\begin{defi} \label{defi-angle} Let $s \in \widetilde{Sp}(H_{{\mathbb
R}},\theta)$. 
\begin{description}
\item[(a)] We will say that $s$ has a non-positive displacement angle
if $s \in \widetilde{Sp}(-)$.
\item[(b)] If $s$ has a non-positive displacement angle define the
displacement angle of $s$ to be the number
\[
\ang(s) := \frac{\op{tr}(\gamma_{x})}{2g} \in (-2\pi,0],
\]
where $x \in \liesp(-2\pi,0]$ is the unique element satisfying
$\op{Exp}(x) = s$.
\item[(c)] For an element  $t \in B_{2g+2}$ of the braid group we will
say that $t$ has a non-positive displacement angle (respectively has
displacement angle $\ang(t)$ equal to $\phi$) if its
image $\tilde{\sigma}_{g}(t) \in
\widetilde{Sp}(H,\theta)$ is contained in the subsemigroup
$\widetilde{Sp}(-)$ (respectively if $\ang(\tilde{\sigma}_{g}(t)) = \phi$).
\end{description}
\end{defi}

\

\begin{rem} \label{rem-compare-genus1}
The notion of a displacement angle that we have just
introduced specializes to the one considered in \cite{abkp} for the
case when $g = 1$. Similarly to the genus one case, the displacement
angle of an element $t \in B_{2g+2}$ makes sense `on the nose'
only if the image of $t$ in $\widetilde{Sp}(H,\theta)$ is contained in 
$\widetilde{U}(g)$. For arbitrary elements $t$ we can talk only about
the direction or the amplitude of a displacement via $t$, but the
actual value of the displacement depends on the point in the Lagrangian
Grassmanian on which $t$ acts. Note also that in contrast with the
$g=1$ case the displacement angle can not be defined directly for
hyperelliptic mapping classes $\tau \in \Delta_{g}^{1}$ but only for
their lifts in $B_{2g +2}$. In other words, rather than working
directly with $\tau$ we need to chose an
element  $t \in B_{2g +2}$ such that $\kappa_{g}^{1}(t) = \tau$ and
work with $t$ instead. One does not see the necessity of such a choice
in the genus one case where the natural map $\kappa_{1,1} : B_{4} =
\widetilde{SL}(2,{\mathbb Z}) \to \Delta_{1,1}$ is an isomorphism and
so we have canonical lifts for Dehn twists. 
\end{rem}

\

\begin{rem} \label{rem-choices} In order to define the subsemigroup
$\widetilde{Sp}(-)$ and to characterize the elements in
$\widetilde{K} \cap \widetilde{Sp}(-)$ in terms of the character
$\widetilde{\det}$ in Lemma~\ref{lem-independence}, we had to make some
(rather mild) choices. Namely we had to choose a maximal compact
subgroup $K \cong U(g) \subset Sp(H_{1}(C_{g},{\mathbb R}))$ and a
$U(g)$-equivariant  
surjection $U(g) \to \Lambda(H_{{\mathbb R}},\theta)$. It is
instructive to examine the geometric meaning of these choices. 

As explained at the beginning of the section, the choice of $K$ is
equivalent to choosing a $\theta$-tamed complex structure $I$ on the real
vector space $H_{{\mathbb R}} := H_{1}(C_{g},{\mathbb R})$. A natural
choice for $I$ will be the Hodge $*$ operator corresponding to a
(conformal class of a) Riemannian metric on $C$.  In other words,
every choice of a complex structure on $C$ corresponds to a choice of
$U(g)$. However not every $I$ comes from a choice of a complex
structure on $C_{g}$. Indeed, specifying the complex structure $I$ on
$H_{{\mathbb R}}$ is equivalent to specifying a splitting of the
complex vector space 
$H_{{\mathbb C}} := H_{{\mathbb R}}\otimes {\mathbb C}$ as
$H_{{\mathbb C}} = H^{{\mathbb C}}\oplus \overline{H^{{\mathbb C}}}$,
i.e. to viewing the triple $(H_{{\mathbb Z}},H_{{\mathbb C}},\theta)$
as a pure polarized Hodge structure of weight one. Thus the choice of
$U(g)$ is equivalent to endowing the torus $H_{{\mathbb
R}}/H_{{\mathbb Z}}$ with the structure of a principally polarized
abelian variety and the choice of $K$ will correspond to a complex
structure on $C_{g}$ if and only if the corresponding period matrix
satisfies the Schottky relations.  

The ambiguity in the choice of the surjection $U(g) \to
\Lambda(H_{{\mathbb R}},\theta)$ also has a transparent geometric
meaning. In order to map $U(g)$ equivariantly to $\Lambda(H_{{\mathbb
R}},\theta)$ we only need to choose a base point $\lambda_{0} \in
\Lambda(H_{{\mathbb R}},\theta)$, i.e. a Lagrangian subspace in 
$H_{{\mathbb R}}$. This choice can be rigidified somewhat if we choose
a Lagrangian subspace in 
$H_{{\mathbb R}}$ which is defined over ${\mathbb Z}$. A standard way
to make such a choice will be to choose a collection of  $a$ and $b$ cycles
on $C_{g}$ and then take $\lambda_{0} = \op{Span}(a_{1}, \ldots, a_{g})$. 
\end{rem}

\

\bigskip

We conclude this section with an estimate for the amplitude of the
displacement angle of some special elements of 
$B_{2g+2}$ which act with fixed points on the
Lagrangian Grassmanian $\Lambda(H_{{\mathbb R}},\theta)$:

\begin{lem} \label{lem-lagrangian-displacement} Let $t \in B_{2g+2}$
be such that $\kappa_{g}^{1}(t) \in \Delta_{g}^{1}$ is a product of
commuting right Dehn twists.  Then $-\pi \leq \ang(t) \leq 0$.
\end{lem}
{\bf Proof.} By hypothesis the element $\sigma_{g}(t)$ is unipotent
and so the element $x \in \liesp(-2\pi,0]$ for which $\op{Exp}(x) =
\tilde{\sigma}_{g}(t)$ must be nilpotent. But the nilpotent elements
in $\liesp(-)$ are contained in the boundary $\partial \liesp(-)$
which is in turn contained in the subdomain 
\[
\liesp[\pi,0] := 
\left\{ x \in \liesp(-) | \op{spectrum}(\gamma_{x}) \subset [-\pi,0] \right\}.
\]
The lemma is proven. \hfill $\Box$

\

\medskip

\begin{rem} \label{rem-fixed} Note that the hypothesis on $t$ in the
previous lemma implies that $\sigma_{g}(t)$ has a fixed
point on $\Lambda(H_{{\mathbb 
R}},\theta)$. In fact by using the description of the boundary
$\partial \liesp(-)$ in Remark~\ref{rem-nilpotent} (ii) one can check that
Lemma \ref{lem-lagrangian-displacement} holds for any non-positive 
element $t$ for
which $\sigma_{g}(t)$ has a fixed
point on $\Lambda(H_{{\mathbb
R}},\theta)$.
\end{rem}

\section{The proof of Theorem~\ref{thm-main}} \label{sec-proof}

Let $f : X \to S^{2}$ be a semistable hyperelliptic symplectic
fibration with a topological section and general fiber $C_{g}$. 

Consider the surjective homomorphisms $\kappa_{g}^{1} : B_{2g + 2} \to
\Delta_{g}^{1}$ and $\sigma_{g} : B_{2g+2} \to Sp(H,\theta)$
introduced in the previous section and let 
\[
{\mathbb K} := \ker\left[B_{2g + 2} \stackrel{\kappa_{g}^{1}}{\to}
\Delta_{g}^{1}\right]. 
\]
Consider the following two elements in $B_{2g+2}$:
\[
h := t_{1}t_{2}\ldots t_{2g+1} \quad \text{and} \quad \bar{h} := 
t_{2g+1}\ldots t_{2}t_{1}.
\]
It is known by \cite{birman-hilden} that the subgroup ${\mathbb K}
\subset B_{2g+2}$ is 
generated by the elements $h^{2g+1}(h\bar{h})^{-1}$ and $h^{2g + 2}$
as a normal subgroup. Furthermore, since $\sigma_{g}({\mathbb K}) =
\{\boldsymbol{1}\}  
\in Sp(H_{{\mathbb
Z}},\theta)$ it follows that 
\[
\tilde{\sigma}_{g|{\mathbb K}} : {\mathbb K} \to {\mathbb Z} \subset 
\widetilde{Sp}(H_{{\mathbb Z}},\theta),
\]
where ${\mathbb Z} = \ker[\widetilde{Sp}(H,\theta) \to
Sp(H,\theta)]$. In particular $\tilde{\sigma}_{g|{\mathbb K}} :
{\mathbb K}
\to {\mathbb Z}$ is a character of ${\mathbb K}$. Our next goal is to
compare the character $\tilde{\sigma}_{g|{\mathbb K}}$ with the character
$\deg_{|{\mathbb K}}$ and the angle character $\ang :
\widetilde{Sp}(-) \to {\mathbb R}$. First we have the following:

\begin{prop} \label{lem-angleofh} \

\begin{itemize}
\item[(i)] The homomorphism $\sigma_{g} : B_{2g+2} \to Sp(H,\theta)$
maps $h$ and $\bar{h}$ to elements of order
$2g+2$ in $Sp(H,\theta)$. In particular $\sigma_{g}(h), \sigma_{g}(\bar{h})
\in U(g)$ for a suitably chosen
maximal compact subgroup $U(g) \subset Sp(H_{{\mathbb R}},\theta)$.
\item[(ii)] The elements $\tilde{\sigma}_{g}(h)$ and
$\tilde{\sigma}_{g}(\bar{h})$ are conjugate in $\widetilde{Sp}(H,\theta)$.
\item[(iii)] The displacement angles of $h$ and $\bar{h}$ are equal to
${\displaystyle \left( - \frac{\pi}{2}\right)}$.
\end{itemize}
\end{prop}
{\bf Proof.} For the proof of parts (i) and (ii) 
we will need the following standard geometric
picture for $h$ and $\bar{h}$. Choose a geometric realization for the
double cover $\nu : C_{g} \to S^{2}$ in which the branch points of
$\nu$ are the $2g+2$ roots of unity $\zeta_{1}, \zeta_{2}, \ldots,
\zeta_{2g+2}$  of order $2g + 2$, labeled consecutively (in the
counterclockwise direction) along the unit
circle.  Let $\op{Diff}^{+}(S^{2},\{\zeta_{i}\}_{i = 1}^{2g + 2})$ be
the group of orientation preserving diffeomorphisms of $S^{2}$ which
leave the set of points $\{\zeta_{i}\}_{i = 1}^{2g + 2}$ invariant and let
$\Gamma_{2g+2} := \pi_{0}(\op{Diff}^{+}(S^{2},\{\zeta_{i}\}_{i =
1}^{2g + 2}))$ be the corresponding mapping class group. It is well known
\cite{birman-hilden,birman-book} that the hyperelliptic mapping
class group $\Delta_{g}$ can be constructed as a central extension
\[
0 \to {\mathbb Z}/2 \to \Delta_{g} \to \Gamma_{2g+2} \to 1,
\]
where the central ${\mathbb Z}/2$ is generated by the mapping class of
the hyper-elliptic involution $\iota$. 
In terms of this realization of $\nu : C_{g} \to S^{2}$ the 
surjective homomorphism 
\[
\xymatrix@1{B_{2g+2} \ar[r] \ar@/_1pc/[rr]_-{\rho_{g}} & \Delta_{g} \ar[r] &
 \Gamma_{2g+2}
}
\]
can be described explicitly
\cite[p.164]{birman-book}: 
\begin{itemize}
\item $\rho_{g} : B_{2g+2} \to \Gamma_{2g+2}$
sends the positive half-twist $t_{i} \in B_{2g+2}$ to
the mapping class $x_{i}$ of a Dehn twist on $S^{2}$ which is the
identity outside of a small neighborhood of the circle segment
connecting $\zeta_{i}$ with $\zeta_{i+1}$  and which switches
$\zeta_{i}$ with $\zeta_{i+1}$. 
\item The kernel $\ker\left[B_{2g+2} \stackrel{\rho_{g}}{\to}
\Gamma_{2g+2}\right]$ is 
generated as a normal subgroup by the elements $h\bar{h}$ and $h^{2g+2}$.
\end{itemize}
Note that $h^{2g+2}$ is a full twist in $B_{2g+2}$ and so in the above
realization of $\nu : C_{g} \to S^{2}$ the mapping class
$\rho_{g}(h^{2g+2})$ can be represented by a rotation on $S^{2}$
through angle $2\pi$. 
In particular we see that $\rho_{g}(h) \in \Gamma_{2g+2}$ is
the mapping class of the counterclockwise ${\displaystyle
\frac{\pi}{g+1}}$-rotation on $S^{2}$. Similarly $\rho_{g}(\bar{h}) \in
\Gamma_{2g+2}$ is
the mapping class of the clockwise ${\displaystyle
\frac{\pi}{g+1}}$-rotation. This proves part (i) of the proposition.

Next observe that $\rho_{g}(h)$ and $\rho_{g}(\bar{h})$
are manifestly conjugate in $\Gamma_{2g+2}$. Indeed we have
$\rho_{g}(\bar{h}) = s\rho_{g}(h)s^{-1}$, where $s$ is the antipodal
involution $s(\zeta) = 1/\zeta$ on $S^{2}$.

The elements $\rho_{g}(h)$ and $s$ generate a dihedral subgroup
$D_{2g+2} \subset \Gamma_{2g+2}$ of $\Gamma_{2g+2}$. The preimage
$\widetilde{D}_{2g+2} := \Delta_{g}\times_{\Gamma_{2g+2}} D_{2g + 2}$ of
$D_{2g + 2}$ in $\Delta_{g}$ is a central extension of $D_{2g+2}$ 
by ${\mathbb Z}/2$, which is the pull back of the standard Heisenberg
extension 
\[
0 \to {\mathbb Z}/2 \to H_{2} \to {\mathbb Z}/2\oplus {\mathbb Z}/2
\to 0,
\]
via the canonical quotient map $D_{2g+2} \twoheadrightarrow  
{\mathbb Z}/2\oplus {\mathbb Z}/2$.

Since the Heisenberg extension splits over any cyclic subgroup, it
follows that we can find a natural lift of $s$ to an element in 
$\widetilde{D}_{2g+2} \subset \Delta_{g}$ which conjugates
$\kappa_{g}(h)$ into $\kappa_{g}(\bar{h})$. Combined with the fact
that $\Delta_{g} \to Sp(H,\theta)$ is a group homomorphism this
implies that $\sigma_{g}(\bar{h}) = s'\sigma_{g}(h)s^{'-1}$ for a
suitably chosen $s' \in Sp(H,\theta)$. Let $\tilde{s} \in
\widetilde{Sp}(H,\theta)$ be an element which maps to 
$s' \in Sp(H,\theta)$. We will check that
$\tilde{s}\tilde{\sigma}_{g}(h)\tilde{s}^{-1} =
\tilde{\sigma}_{g}(\bar{h})$.  First observe that for any symplectic
transvection $t \in Sp(H_{{\mathbb R}},\theta)$ the conjugate element 
$s'ts^{'-1} \in Sp(H_{{\mathbb R}},\theta)$ is also a symplectic
transvection. Let $\tilde{t} \in
\widetilde{Sp}(H_{{\mathbb R}},\theta)$ be the standard lift of $t$
described in Lemma~\ref{lem-lift}. By definition $\tilde{t}$ has fixed points
on $\widetilde{\Lambda}(H_{{\mathbb R}},\theta)$ and so
$\tilde{s}\tilde{t}\tilde{s}^{-1}$ will also have fixed points. Hence
Lemma~\ref{lem-lift} we conclude that
$\tilde{s}\tilde{t}\tilde{s}^{-1}$ is the standard lift of the
transvection $s'ts^{'-1}$. Consider now the element $\sigma_{g}(h) =
\prod_{i = 1}^{2g+1} \sigma_{g}(t_{i})$. By Corollary
\ref{cor-homomorphism} we know that $\tilde{\sigma}_{g}(t_{i}) =
\widetilde{\sigma_{g}(t_{i})}$ is the standard lift of the
transvection $\sigma_{g}(t_{i})$ and that $\tilde{\sigma}_{g}(h) =
\prod_{i = 1}^{2g+1} \tilde{\sigma_{g}}(t_{i})$. In particular we have
\[
\begin{split}
\tilde{s}\tilde{\sigma}_{g}(h)\tilde{s}^{-1} & = \prod_{i = 1}^{2g+1}
(\tilde{s}\tilde{\sigma_{g}}(t_{i}) \tilde{s}^{-1}) = \prod_{i =
1}^{2g+1} \widetilde{s't_{i}s^{'-1}} = \prod_{i =
1}^{2g+1} \tilde{\sigma}_{g}({s't_{i}s^{'-1}}) = \\
& = \tilde{\sigma}_{g}\left(\prod_{i = 1}^{2g+1}(s't_{i}s^{'-1})\right) = 
\tilde{\sigma}_{g}(s'hs^{'-1}) = \tilde{\sigma}_{g}(\bar{h}).
\end{split}
\] 
This completes the proof of part (ii) of the proposition. 

We are now ready  to prove part (iii). The elements $h,
\bar{h} \in \Delta_{g}^{1}$ are products of right handed Dehn twists
and so the elements $\tilde{\sigma}_{g}(h),
\tilde{\sigma}_{g}(\bar{h})$ both belong to
$\widetilde{Sp}(-)$. Since by part (i) of the proposition we know that 
$\tilde{\sigma}_{g}(h),
\tilde{\sigma}_{g}(\bar{h})$ also belong to
$\widetilde{U}(g)$ it follows that $h$ and $\bar{h}$ have well defined
negative displacement angles. Furthermore, note that $h,
\bar{h}$ must
have the same displacement angle since $\tilde{\sigma}_{g}(h)$ and 
$\tilde{\sigma}_{g}(\bar{h})$ belong to the same
conjugacy class in $\widetilde{U}(g)$. In view of this and the fact
that  $\widetilde{\det}$ is additive on $\widetilde{U}(g)$, it
suffices to show that the displacement angle of $h\bar{h}$ is equal to
$-\pi$. Since the element $h\bar{h}$ maps to the mapping class in
$\Delta_{g}^{1}$ represented by the hyperelliptic involution $\iota :
C_{g} \to C_{g}$, it follows that $\sigma_{g}(h\bar{h}) = -
\boldsymbol{1} \in Sp(H,\theta) \subset SL(H)$. In particular
$\tilde{\sigma}_{g}(h\bar{h}) \in \widetilde{Sp}(H,\theta)$ will be
the unique lift of  $-\boldsymbol{1}$,  which has fixed  points on
$\widetilde{\Lambda}(H_{{\mathbb R}},\theta)$. But we already know one
such lift of $-\boldsymbol{1}$, namely the negative generator ${\mathfrak c}$
of the center  of $\widetilde{Sp}(H,\theta)$. Indeed the center
$Z(\widetilde{Sp}(H,\theta))$ of $\widetilde{Sp}(H,\theta)$ is an
infinite cyclic group which maps onto
the center $Z(Sp(H,\theta)) \cong {\mathbb Z}/2$ of
$Sp(H,\theta)$. Since the latter is generated by $-\boldsymbol{1}$ we
conclude that $\tilde{\sigma}_{g}(h\bar{h}) = {\mathfrak c}$. Finally
the element  $-\boldsymbol{1}$ considered as an element in $U(g)$ has
eigenvalues $e^{\pi i}$ and so $\widetilde{\det}({\mathfrak c}) = -
g\pi$. The proposition is proven. \hfill $\Box$

\bigskip

We can now finish the proof of Theorem~\ref{thm-main}. Let ${\mathfrak
t}_{1}\ldots {\mathfrak t}_{N} = 1$ be the relation in
$\Delta_{g}^{1}$ corresponding to the monodromy representation of the
pencil $f : X \to S^{2}$. Choose elements $\nu_{1}, \ldots, \nu_{N} \in
B_{2g+2}$ satisfying $\kappa_{g}^{1}({\mathfrak t}_{i}) =
\nu_{i}$. Then the element $\mu := \nu_{1}\ldots \nu_{N} \in B_{2g+2}$
belongs to the subgroup ${\mathbb K} \subset B_{2g+2}$. But the
subgroup ${\mathbb K}$ is normally generated by the elements
$h^{2g+1}(h\bar{h})^{-1}$ and $h^{2g+2}$ and hence is contained in the
subgroup ${\mathbb L}$ of $B_{2g+2}$ which is normally generated by
the elements $h$ and $\bar{h}$. Let now ${\mathbb L}(-) = {\mathbb
L}\cap \tilde{\sigma}_{g}^{-1}(\widetilde{Sp}(-))$ be the non-positive
sub semigroup of ${\mathbb L}$. Then we have two natural ${\mathbb
R}_{\leq 0}$-valued characters on ${\mathbb L}(-)$:
\[
-\deg : {\mathbb L}(-) \to {\mathbb Z}_{\leq 0} \subset {\mathbb
R}_{\leq 0} \quad \text{and} \quad \chi := \ang\circ
\tilde{\sigma}_{g} : {\mathbb L}(-) \to 
{\mathbb R}_{\leq 0}.
\]
By the previous proposition we have $\chi(h) =
\chi(\bar{h}) = - \pi/2$. On the other hand $\deg(h) = \deg(\bar{h}) =
2g+1$  and since ${\mathbb L}$ is normally generated by $h$ and
$\bar{h}$ we have 
\[
\chi = - \frac{\pi}{4g+2}\cdot \deg
\]
on ${\mathbb L}(-)$.

Since $\mu \in {\mathbb L}(-)$ this implies 
\[
\ang(s_{1}s_{2}\ldots s_{N}) = -\frac{\pi N}{4g+2},
\]
where $s_{i} \in \widetilde{Sp}(-)$ is the canonical lift (see Lemma
\ref{lem-lift}) of the simple homological Dehn twist
$\sigma_{g}({\mathfrak t}_{i})$. Finally, the product $\prod_{i =
1}^{N} s_{i}$ can be rewritten as $\prod_{j = 1}^{D} p_{j}$ where
$p_{j}$ is the product of those Dehn twists among the $s_{i}$'s which
correspond to the all the cycles vanishing at the $j$-th singular
fiber of $f : X \to S^{2}$. Since each $p_{j}$ is a product of
commuting Dehn twists Lemma \ref{lem-lagrangian-displacement} applies
and so $\ang(p_{j}) \geq -\pi$. Thus we get 
\[
\ang\left(\prod_{j = 1}^{D} p_{j}\right) = \sum_{j =1}^{D} \ang(p_{j})
\geq -\pi, 
\]
and so
\[
-D\pi \leq - \frac{\pi N}{4g +2}.
\]
Theorem~\ref{thm-main} is proven. \hfill $\Box$

\section{Concluding remarks} \label{sec-remarks}

The  method of proof of Theorem~\ref{thm-main} can be generalized in several
directions and we intend to pursue such generalizations in a forthcoming
paper. We conclude the present discussion by indicating some of
the possible venues of generalization:

\

\bigskip

\noindent
{\bf (i)}  Obtain analogues of the Szpiro inequality for a general
semistable family of curves of genus $g$. 

The inequality which our method produces in this case is similar to
the one in Theorem~\ref{thm-main} but instead
of $N$ - the length of the word of the Dehn twists defining the 
fibration we have additional characteristics   - the number of lantern
relations in the monodromy word. As of now the only missing element is
the exact coefficient for the number of lantern relations in the
word.

In order to explain this in more detail recall first 
that the main obstruction to generalizing the genus one proof
from \cite{abkp} to higher genus is the existence of two different natural 
central extensions of $\op{Map}_{g,1}$ when $g \geq 2$. This leads to
two different and rather mild geometric constraints on the monodromy
representation of a Lefschetz family, which when taken separately are
not stringent enough to produce an inequality of Szpiro type.

The idea therefore is to look for a different group which surjects
onto $\op{Map}_{g,1}$ and for which the pullback of the two central
extensions  of $\op{Map}_{g,1}$ become easier to compare,
i.e. coincide or at least become proportional. The fact that 
in the hyperelliptic case one can successfully implement this idea  by
using the braid group suggests that for a general family of genus $g$
curves one may try to work with the Artin braid group
$\op{Art}_{\Gamma}$ corresponding to the $T$-shaped graph $\Gamma$ of
Wajnryb's presentation of $\op{Map}_{g}^{1}$ 
\cite{wajnryb}. Since $\Gamma$ is a tree we can
again apply Lemma~\ref{lem-vanish} and conclude that 
pullback of the central extension \eqref{eq-centralR}
will be zero in $H^2(\op{Art}_{\Gamma},{\mathbb Z})$.

Next observe that in Wajnryb's presentation of $\op{Map}_{g}^{1}$ an
extra relation - the {\em lantern raltion} - appears  when $g\geq  3$.
The kernel ${\mathbb K}_{\Gamma} := \ker[\op{Art}_{\Gamma} \to
\op{Map}_{g}^{1}]$  is now normally generated by the 
hyperelliptic relations $h^{2g+1}(h\bar{h})^{-1}$ and $h^{2g + 2}$ in a subgroup $B_{2g+2} \subset \op{Art}_{\Gamma}$
and by the lantern relation. In this setup one can argue that  
the analogue ${\mathbb K}_{\Gamma} \to
{\mathbb Z}$ of the character $\tilde{\sigma}_{g|{\mathbb K}}$ can be
expressed as $a \deg + b \ell$ where $\ell$  is a character which is 
non-trivial on the lantern relations in ${\mathbb K}_{\Gamma}$.

This reasoning leads to a formula for the left hand side of the Szpiro
inequality  which instead  of one parameter $N$ involves two parameters
$N$ and  $L$ where the latter corresponds to the number
of lantern relations in the monodromy word in $\op{Map}_{g}^{1}$.
We do not know the exact coefficient for $L$ yet but it is clear that
this coefficient can not be trivial. Indeed otherwise the
character ${\mathbb K}_{\Gamma} \to {\mathbb Z}$ will be (rationally)
proportional  $\deg$ which contradicts the fact that the central
extension \eqref{eq-centralR}  pull back to a non-torsion element in
$H^{2}(\op{Map}_{g}^{1},{\mathbb Z})$.
 
\bigskip

\noindent
{\bf (ii)}  Work out the hyperelliptic Szpiro inequality in 
the case of a general base curve. 

In this setup one needs to change the right hand
side of the inequality  by an appropriate multiple of  the genus
of the base curve (in the case of elliptic fibrations the correct 
modification was worked out by S.-W. Zhang \cite{zhang}). 

To carry out the argument in this case  one needs to analyze 
more carefully the monodromy word corresponding to a Lefschetz  family of
surfaces over a closed base surface $C$ of arbitrary genus. The product of
local monodromy  transformations ${\mathfrak t}_i$ is not $1$ but a
product of at most $g(C)$ commutators $aba^{-1}b^{-1}$. Now the set of
commutators in $\widetilde{Sp}(2g,{\mathbb R})$ maps isomorphically
to the corresponding set in $Sp(2g,{\mathbb R})$ and in fact
consists of the elements $h\in \widetilde{Sp}(2g,{\mathbb R})$ with
bounded displacement angle. Thus there is a number $A$ so that
for any $\lambda$ in the Lagrangian Grassmanian, we have
that the angle between $h\cdot \lambda$ and $\lambda$ is not less than
$-A\pi$. This results in a correction term $Ag(C)$ in the right hand
side of the Szpiro formula.

\

\bigskip

\noindent
{\bf (iii)} Find a monodromy proof of the negativity of
self intersections for sections in a symplectic Lefschetz fibration. 

The existing proofs of this fact use either the theory of
pseudo-holomorphic curves or properties of the Seiberg-Witten Floer
homology.  Both  apparoaches are quite technical and require the existence 
of global solutions of special partial differential equations.

In contrast the proof suggested by our method is completely algebraic
and relies 
only on the combinatorial properties of the monodromy group. The idea
is to look for a `displacement angle' 
 description of the subsemigroup in 
$\op{Map}_{g}^{1}$  consisting of elements which are negative  with respect to
the Morita move on the boundary circle for the uniformizing disk of
the fiber (see e.g. \cite{abkp}). Here is a brief sketch of such a proof.

Let $\rho : \op{Map}_{g}^{1} \to
\op{Homeo}^{+}(S^{1})$ be the Morita homomorphism
and let $e \in H^{2}(\op{Map}_{g}^{1},{\mathbb Z})$ be the pullback of
the natural central extension
\[
0 \to {\mathbb Z} \to \op{Homeo}^{\op{per}}({\mathbb R})
\to \op{Homeo}^{+}(S^1) \to 1
\]
via $\rho$. Recall that the middle term
of the extension $e$
can be naturally identified with the mapping class group
$\op{Map}_{g,1}$ \cite{morita-characteristic}.

Using Lemma \ref{lem-vanish} one can again argue that the pullback of
the extension $e$ to $\op{Art}_{\Gamma}$ is trivial and so there is
a  natural homomorphism $\op{Art}_{\Gamma} \to
\op{Homeo}^{\op{per}}({\mathbb R})$. Now by shearing the
hyperbolic disk one can show (see Proposition 2.1 and remark 2.2 in
\cite{smith}) that the generators in $\op{Art}_{\Gamma}$ covering
right handed Dehn twists in $\op{Map}_{g}^{1}$ get mapped to elements 
in $\op{Homeo}^{\op{per}}({\mathbb R})$ which fix a countable set of
points in ${\mathbb R}$ and move all other points in ${\mathbb R}$ to
the right. In particular any positive word $x$ in the generators of
$\op{Art}_{\Gamma}$ which is in $\ker[\op{Art}_{\Gamma} \to
\op{Map}_{g}^{1}]$ defines a nontrvial right-shifting periodic homeomorphism  
of ${\mathbb R}$.

Given a symplectic Lefschetz family over $S^{2}$ with a fixed
symplectic section $s$, we get a lifting of the
monodromy representation into the group $\op{Map}_{g}^{1}$ and
the selfintersetion of $s$ corresponds to the image of (the lifting in
$\op{Art}_{\Gamma}$ of) the monodromy word in
$\op{Homeo}^{\op{per}}({\mathbb R})$. Thus the fact that the resulting
element is right-shifting corresponds precisely 
to the negativity of the selfintersection of $s$.

We ca make this more precise (in the spirit of Lemma
\ref{lem-lagrangian-displacement}) if we manage  to bound the amount
of the shifting in terms of the numbers $N$ and $L$ mentioned in (i)
above. This will give an explicit bound on the selfintersection
of $s$ in terms of $N$ and $L$. Some preliminary computations we have
made show that the $N$-$L$ bounds one gets imply also effective bounds 
on $s^{2}$ which depend on $N$ only. Estimates of this type are of
independent interest since they provide a simple way
to show finiteness of types of symplectic Lefshetz pencils ($s^2 =
-1$) for a given genus $g > 1$.

{\bf (iv)} Find arithmetic analogues of the symplectic Szpiro
inequality.

It is very tempting to try to apply our method  to
the arithmetic situation. 

In this  case the analogue of the  global displacement angle is
clearly the height  of a curve with a point. The monodromy representation
corresponds to the Galois representation.
There are several problems with this approach.
One of them is that there are no direct analogs of the
${\mathbb Z}$-central extension in this case and the other one
is that we do not  yet understand what should the local inequalities be
near singular fibers.

\end{document}